\documentclass{article}

\usepackage{amsfonts}
\usepackage{amsmath}
\usepackage{amssymb}
\usepackage{latexsym}

\newtheorem{theorem}{Theorem}[section]
\newtheorem{lemma}[theorem]{Lemma}

\newtheorem{proposition}[theorem]{Proposition}

\newtheorem{corollary}[theorem]{Corollary}

\newtheorem{remark}[theorem]{Remark}

\numberwithin{equation}{section}

\newcommand{\R}{\mathbb{R}}

\begin{document}

\title
{The Index of a Vector Field on an Orbifold with Boundary}

\author{Elliot Paquette, Kalamazoo College
\thanks{The first author was supported by a Kalamazoo College Field Experience grant.} \\
1880 Sumac Ave.               \\
Boulder, CO 80304               \\
elliot.paquette@gmail.com     \\\\
Christopher Seaton
\thanks{The second author was supported by a Rhodes College Faculty
Development Endowment Grant.}            \\
Department of Mathematics and Computer Science \\
Rhodes College \\2000 N. Parkway       \\ Memphis, TN 38112    \\
seatonc@rhodes.edu }

\maketitle

Subject class: Primary 57R25, 57R12, 55R91


\begin{abstract}

A Poincar\'{e}-Hopf theorem in the spirit of Pugh is proven for
compact orbifolds with boundary.  The theorem relates the index sum
of a smooth vector field in generic contact with the boundary
orbifold to the Euler-Satake characteristic of the orbifold and a
boundary term. The boundary term is expressed as a sum of Euler
characteristics of tangency and exit-region orbifolds.  As a
corollary, we express the index sum of the vector field induced on
the inertia orbifold to the Euler characteristics of the associated
underlying topological spaces.

\end{abstract}


\section{Introduction}
\label{sec-intro}

Pugh gave a generalization (see \cite{Pugh}) of the Poincar\'e-Hopf
Theorem for manifolds with boundary for continuous vector fields in
generic contact with the boundary.  This generalization bears the
elegance of associating the index sum with a sum of Euler
characteristics only.  Here, we will show that, in the case of a
compact orbifold with boundary and a smooth vector field in generic
contact with the boundary, Pugh's result extends naturally. A proper
introduction to orbifolds and the precise definition we will use are
available as an appendix in \cite{chenruangwt}.  Note that this
definition of an orbifold requires group actions to have fixed-point
sets of codimension 2 as opposed to, e.g., \cite{thurston}; we make
this requirement as well. By ``smooth," we will always mean
${\mathcal C}^\infty$.

The main result we will prove is as follows.

\begin{theorem}
\label{thrm-mainresult}

Let $Q$ be an $n$-dimensional smooth, compact orbifold with
boundary.  Let $Y$ be a smooth vector field on $Q$ that is in
\emph{generic contact} with $\partial Q$, and then
\begin{equation}
\label{eq-mainformula}
    \mathfrak{Ind}^{orb}(Y; Q) = \chi_{orb}(Q,
    \partial Q) + \sum_{i = 1}^n \chi_{orb}(R_{-}^i, \Gamma^i).
\end{equation}

\end{theorem}

The expressions $\mathfrak{Ind}^{orb}$ and $\chi_{orb}$ are the
orbifold analogues of the manifold notions of the topological index
of a vector field and the Euler characteristic, respectively. The
definitions of both of these, along with the orbifolds $R_-^i$,
$\Gamma^i$, and generic contact, will be reviewed in Section
\ref{sec-preliminaries}.

In this paper, we will follow a procedure resembling Pugh's original
technique, and we will show that many of the same techniques
applicable to manifolds can be applied to orbifolds as well.  In
Section \ref{sec-preliminaries} we will explain our notation and
review the result of Satake's which relates the orbifold index to
the Euler-Satake characteristic for closed orbifolds.  We give the
definition of each of these terms.  In Section \ref{sec-double}, we
will show that a neighborhood of the boundary of an orbifold may be
decomposed as a product $\partial Q \times [0, \epsilon)$. We will
then construct the double of $Q$ and charts near the boundary
respecting this product structure.  This will generalize well-known
results and constructions for manifolds with boundary.
 Section \ref{sec-morseindex} provides elementary
results relating the topological index of an orbifold vector field
to an orbifold Morse Index.  The orbifold Morse Index is defined in
terms of the Morse Index on a manifold in a manner analogous
Satake's definition of the topological vector field index.  These
results generalize corresponding results for manifolds.  In Section
\ref{sec-mainresult}, we will use the above constructions to show
that a smooth vector field on $Q$ may perturbed near the boundary to
form a smooth vector field on the the double whose index can be
computed in terms of the data given by the original vector field. We
use this to prove Theorem \ref{thrm-mainresult}.  We also prove
Corollary \ref{cor-inertiaversion}, which gives a similar formula
where the left-hand side is the orbifold index of the induced vector
field on the inertia orbifold and on the right-hand side, the
Euler-Satake characteristics are replaced with the Euler
characteristics of the underlying topological spaces.

Another generalization of the Poincar\'{e}-Hopf Theorem to orbifolds
with boundary follows as a corollary to Satake's Gauss-Bonnet
Theorem for orbifold with boundary in \cite{satake2}; this and
related results are explored in \cite{seaton1}. In each of these
cases, the boundary term is expressed by evaluation of an auxiliary
differential form representing a global topological invariant of the
boundary pulled back via the vector field.  The generalization given
here expresses the boundary term in terms of Euler-Satake
characteristics of suborbifolds determined by the vector field.

The first author would like to thank Michele Intermont for guiding
him through much of the background material required for this work.
The second author would like to thank Carla Farsi for helpful
conversations and suggesting this problem.


\section{Preliminaries and Definitions}
\label{sec-preliminaries}

Satake proved a Poincar\'{e}-Hopf Theorem for closed orbifolds;
however, he worked with a slightly different definition of orbifold,
the so-called $V$-manifold (see \cite{satake1} and \cite{satake2}).
A $V$-manifold corresponds to modern day \emph{effective} or
\emph{reduced} (codimension-2) orbifold. An effective orbifold is
such that the group in each chart acts effectively (see
\cite{chenruangwt}). We will adapt the language of his result and
use it here.

\begin{theorem}[Satake's Poincar\'{e}-Hopf Theorem for Closed Orbifolds]
\label{thrm-satakeindex}

Let $Q$ be an effective, closed orbifold, and let $X$ be a vector
field on $Q$ that has isolated zeroes.  Then the following
relationship holds.
\[
    \mathfrak{Ind}^{orb}(X; Q) =
    \chi_{orb}(Q)
\]

\end{theorem}

Note that the requirement that $Q$ is effective is unnecessary; as
mentioned in \cite{chenruangwt}, an ineffective orbifold can be
replaced with an effective orbifold $Q_{red}$, and the differential
geometry of the tangent bundle (or any other \emph{good} orbifold
vector bundle) is unchanged.

The \emph{orbifold index} $\mathfrak{Ind}^{orb}(X; p)$ at a zero $p$
of the vector field $X$ is defined in terms of the topological index
of a vector field on a manifold. Let a neighborhood of $p$ be
uniformized by the chart $\{V, G, \pi\}$ and choose $x \in V$ with
$\pi(x) = p$. Let $G_x \leq G$ denote the isotropy group of $x$.
Then $\pi^\ast X$ is a $G$-equivariant vector field on $V$ with a
zero at $x$.  The orbifold index at $p$ is then defined as
\[
    \mathfrak{Ind}^{orb}(X; p)
    =
    \frac{1}{|G_x|}\mathfrak{Ind}\left(\pi^\ast X;
    x\right),
\]
where $\mathfrak{Ind}\left(\pi^\ast X; x\right)$ is the usual
topological index of the vector field $\pi^\ast X$ on the manifold
$V$ at $x$ (see \cite{GuilleminPollack} or \cite{Milnor}). Note that
this definition does not depend on the chart, nor on the choice of
$x$. We use the notation
\[
    \mathfrak{Ind}^{orb}(X; Q) = \sum\limits_{p \in Q, X(p) = 0} \mathfrak{Ind}^{orb} (X;
    p).
\]

The \emph{Euler-Satake characteristic} $\chi_{orb}(Q)$ is most
easily defined in terms of an appropriate simplicial decomposition
of $Q$. In particular, let ${\mathcal T}$ be a simplicial
decomposition of $Q$ so that that the isomorphism class of the
isotropy group is constant on the interior of each simplex (such a
simplicial decomposition always exists; see \cite{moerdijkpronk}).
For the simplex $\sigma$, the (isomorphism class of the) isotropy
group on the interior of $\sigma$ will be denoted $G_\sigma$. The
Euler-Satake characteristic of $Q$ is
\[
    \chi_{orb}(Q)
    = \sum_{\sigma \in {\mathcal T}}(-1)^{\mbox{\scriptsize dim\:} \sigma}\frac{1}{|G_\sigma|}.
\]
This coincides with Satake's \emph{Euler characteristic of $Q$ as a
$V$-manifold}.  Note that it follows from this definition that if $Q
= Q_1 \cup Q_2$ for orbifolds $Q_1$ and $Q_2$ with $Q_1 \cap Q_2$ a
suborbifold, then
\begin{equation}
\label{eq-additiveeulerchar}
    \chi_{orb}(Q) = \chi_{orb}(Q_1) + \chi_{orb}(Q_2) - \chi_{orb}(Q_1 \cap Q_2).
\end{equation}
In the case that $Q$ has boundary, $\chi_{orb}(Q)$ is defined in the
same way.  We let
\[
    \chi_{orb} (Q, \partial Q) = \chi_{orb}(Q) - \chi_{orb}(\partial Q).
\]
This coincides with Satake's \emph{inner Euler characteristic of $Q$
as a $V$-manifold with boundaries}.  The reader is warned that there
are many different Euler characteristics defined for orbifolds; both
the topological index of a vector field and Euler-Satake
characteristic used hear are generally rational numbers.

Vector fields in \emph{generic contact} have orbifold exit regions,
which we will now describe.  Let $Q$ be a compact $n$-dimensional
orbifold with boundary.  In Section \ref{sec-double} Lemma
\ref{lm-normalspace}, we will show that, as with the case of
manifolds, there is a neighborhood of $\partial Q$ in $Q$
diffeomorphic to $\partial Q \times [0, \epsilon)$.  Given a metric,
he tangent bundle of $Q$ on the boundary decomposes with respect to
this product so that there is a well-defined normal direction at the
boundary. Let $R_{-}^1$ be the closure of the subset of $\partial Q$
where the vector field points out of $Q$. Analogously, let $R_{+}^1$
be the closure of the subset of $\partial Q$ where the vector field
points into $Q$.  We require that $R_-^1$ and $R_+^1$ are orbifolds
with boundary of dimension $n - 1$.  The subset of $\partial Q$
where the vector field is tangent to $\partial Q$ is denoted
$\Gamma^1$; we require that $\Gamma^1$ is a suborbifold of $\partial
Q$ of codimension $1$. Note that, by the continuity of the vector
field, the component of the vector field pointing outward must
approach zero near the boundary of $R_{-}^i$ and $R_{+}^i$. Hence
$\Gamma^1 = \partial R_{-}^1 = \partial R_{+}^1$.

The vector field is tangent to $\Gamma^1$, and so it may be
considered a vector field on the orbifold $\Gamma^1$.  We again
require this vector field to have orbifold exit regions. Call
$R_{-}^2$ the closure of the subset of $\partial \Gamma^1$ where the
vector field points out of $R_-^1$, and $R_{+}^2$ the closure of the
subset where it points into $R_-^1$. The subset of $\Gamma^1$ where
the vector field is tangent to $\Gamma^1$ is denoted $\Gamma^2$, and
is required to be a codimension-$1$ suborbifold of $\Gamma^1$.

In the same way, we define $\Gamma^i$, $R_{-}^i$, $R_{+}^i$,
requiring that these sets form a chain of closed suborbifolds
$\{\Gamma^i\}_{i=1}^n$ and compact orbifolds with boundary
$\{R_{-}^i\}_{i=1}^n$.  We require that that $\mbox{dim}\: R_{-}^i =
\mbox{dim}\: R_{+}^i =  n - i$ and $\mbox{dim}\:\Gamma^i = n - i -
1$.  Since each successive $\Gamma^i$ will have strictly smaller
dimension, we eventually run out of space, and so both of these
sequences terminate.  The last entry in the sequence of $\Gamma^i$
will be $\Gamma^n$, which is necessarily the empty set.


\section{Formation of the Double Orbifold}
\label{sec-double}

In the proof of Theorem \ref{thrm-mainresult}, we will pass from an
orbifold with boundary to a closed orbifold in order to employ
Theorem \ref{thrm-satakeindex}.  In this section we will construct
the double of an orbifold with boundary.  In the process, we will
develop charts near the boundary of a specific form which will be
required in the sequel. The construction of the double is similar to
the case of a manifold; see \cite{munkresdifftop}.

Let $\mathbf{B}_x(r)$ denote the ball of radius $r$ about $x$ in
$\R^n$ where $\R^n$ has basis $\{ e_i \}_{i=1}^n$.  For convenience,
$\mathbf{B}_0$ will denote the ball of radius $1$ centered at the
origin in $\R^n$. We let $\R_+^n = \{ x_1, \ldots , x_n : x_n \geq 0
\}$ where the $x_i$ are the coordinates with respect to the basis
$\{ e_i \}$, $\mathbf{B}_x^+(r) = \mathbf{B}_x(r) \cap \R_+^n$, and
$\mathbf{B}_0^+ = \mathbf{B}_0 \cap \R_+^n$. Also $\mathbf{B}_0^k$
will denote the ball of radius $1$ about the origin in $\R^k$.

Let $Q$ be a compact orbifold with boundary.  For each point $p \in
Q$, we choose an orbifold chart $\{ V_p, G_p, \pi_p \}$ where $V_p$
is $\mathbf{B}_0$ or $\mathbf{B}_0^+$ and $\pi_p(0) = p$. Let $U_p$
denote $\pi_p(V_p) \subseteq Q$ for each $p$, and then the $U_p$
form an open cover of $Q$. Choose a finite subcover of the $U_p$,
and on each $V_p$ corresponding to a $U_p$ in the subcover, we put
the standard Riemannian structure on $V_p$ so that the
$\left\{\frac{\partial}{\partial x_i} \right\}$ form an orthonormal
basis.  Endow $Q$ with a Riemannian structure by patching these
Riemannian metrics together using a partition of unity subordinate
to the finite subcover of $Q$ chosen above.

Now, let $p \in Q$, and then there is a geodesic neighborhood $U_p$
about $p$ uniformized by $\{ V_p, G_p, \pi_p \}$ where $V_p =
\mathbf{B}_0(r)$ or $\mathbf{B}_0^+(r)$ for some $r
> 0$ where $G_p$ acts as a subgroup of $O(n)$ (see
\cite{chenruangwt}). Identifying $V_p$ with a subset of $T_0 V_p$
via the exponential map, we can assume as above that $\{ e_i \}$
forms an orthonormal basis with respect to which coordinates will be
denoted $\{x_i \}$. In the case with boundary, $\mathbf{B}_0^+(r)$
corresponds to points with $x_n \geq 0$.  We call such a chart a
\emph{geodesic chart of radius $r$ at $p$}.  Note that in such
charts, the action of $\gamma \in G_p$ on $V_p$ and the action of
$d\gamma = D(\gamma)_0$ on a neighborhood of $0$ in $T_0 V_p$ (or in
half-space in the case with boundary) are identified via the
exponential map.

We may now introduce the following lemma.

\begin{lemma}
\label{lm-normalspace} At every point $p$ in $\partial Q$, there a
geodesic chart at $p$ of the form $\{V_p, G_p, \pi_p\}$ where $G_p$
fixes $e_n$. On the boundary, the tangent space $TQ|_{\partial Q}$
is decomposed orthogonally into $(T\partial Q) \oplus \nu$ where
$\nu$ is a trivial $1$-bundle on which each group acts trivially.

\end{lemma}

{\it \noindent Proof:} Let $p \in \partial Q$, and let a
neighborhood of $p$ be uniformized by the geodesic chart $\{ V_p,
G_p, \pi_p \}$ so that $V_p = \mathbf{B}_0^+(r)$.  Let $\langle
\cdot , \cdot \rangle_0$ denote the inner product on $T_0 V_p$.  Let
$T_0^+$ correspond to the half-space $T_0^+$ in $T_0 V_p$
corresponding to vectors with non-negative $\frac{\partial}{\partial
x_n}$-component.  The exponential map identifies an open ball about
$0 \in T_0^+$ with $V_p$.

Suppose that $\gamma$ is an arbitrary element of $G_p$ so that
$d\gamma$ acts on $T_0 V_p$. Any $v \in T_0^+$ satisfies
$\left\langle v, \frac{\partial}{\partial x_n} \right\rangle_0 \geq
0$. Furthermore, $(d\gamma) v \in T_p^+$, so $\left\langle (d\gamma)
v, \frac{\partial}{\partial x_n} \right\rangle_0 \geq 0$, or
equivalently, $\left\langle v, d\gamma^{-1} \frac{\partial}{\partial
x_n} \right\rangle_0 \geq 0$ for all $v \in T_0^+$.

It will be shown that this implies $G_p$ fixes
$\frac{\partial}{\partial x_n}$.  Pick $j\neq n$; since
$\frac{\partial}{\partial x_j} \in T_0^+$, $\left\langle
\frac{\partial}{\partial x_j}, d\gamma^{-1} \frac{\partial}{\partial
x_n} \right\rangle_0 \geq 0$. However, $-\frac{\partial}{\partial
x_j}$ is also a vector in $T_0^+$, and so $\left\langle
-\frac{\partial}{\partial x_j}, d\gamma^{-1}\frac{\partial}{\partial
x_n} \right\rangle_0 \geq 0$. By the linearity of the inner product,
this is only possible if $\left\langle \frac{\partial}{\partial
x_j}, d\gamma^{-1}\frac{\partial}{\partial x_n} \right\rangle_0 =
0$. Furthermore, since $j \neq n$ was arbitrary, this implies that
$d\gamma^{-1} \frac{\partial}{\partial x_n}$ has no component in the
direction of any $\frac{\partial}{\partial x_j}$, $j\neq n$. Since
$d\gamma^{-1}$ is an isometry, $d\gamma^{-1}
\frac{\partial}{\partial x_n} = \pm \frac{\partial}{\partial x_n}$,
but because $d\gamma^{-1}T_0^+ = T_0^+$, it must be the case that
$d\gamma^{-1}\frac{\partial}{\partial x_n} =
\frac{\partial}{\partial x_n}$. As $\gamma \in G_p$ was arbitrary,
this implies $G_p$ fixes $\frac{\partial}{\partial x_n}$.

Now, for each $p \in \partial Q$, pick a geodesic chart $\{V_p, G_p,
\pi_p \}$ at $p$, and let $N_p$ denote the constant vector field
$\frac{\partial}{\partial x_n}$ on $V_p$.  Recall from
\cite{satake2} that $\tilde{T}_0V_p$ denotes the $dG_p$-invariant
tangent space of $T_0V_p$ on which the differential of $\pi_p$ is
invertible.  If $q \in \pi_p (V_p) \subset Q$ with geodesic chart
$\{ V_q, G_q, \pi_q \}$ at $q$, then the fact that $D(\pi_q)_p^{-1}
\circ D(\pi_p)_0 : \tilde{T}_0 V_p \rightarrow \tilde{T}_0 V_q$ maps
$\tilde{T}_0 \partial V_p$ to $\tilde{T}_0 \partial V_q$ and
preserves the metric ensures that the value of $N_q(0)$ coincides
with that of $D(\pi_q)_p^{-1} \circ D(\pi_p)_0[N_p(0)]$ up to a
sign. The sign is characterized by the property that for any curve
$c : (-1, 1) \rightarrow V_p$ with derivative $c^\prime(t) = N_p$,
there is an $\epsilon > 0$ such that $c(t)$ is in the interior of
$V_p$ for $t \in (0, \epsilon)$; a curve in $V_q$ with derivative
$D(\pi_q)_p^{-1} \circ D(\pi_p)_0[N_p(0)]$ has the same property.
With this, we see that the $N_p$ patch together to form a
nonvanishing section of $TQ|_{\partial Q}$ that is orthogonal to
$T\partial Q$ at every point; hence, it defines a trivial subbundle
$\nu$ orthogonal to $T\partial Q$. Clearly, $TQ = (T\partial
Q)\oplus \nu$.

{\hfill $\square$}

Let $Q^\prime$ be an identical copy of $Q$.  In order to form a
closed orbifold from the two, the boundaries of these two orbifolds
will be identified via
\begin{equation}
    \label{eq-glue}
    \partial Q \ni x \longleftrightarrow x^\prime \in \partial Q^\prime
\end{equation}
The resulting space inherits the structure of a smooth orbifold from
$Q$ as will be demonstrated below.

First, note that for each point $p \in \partial Q$, by Lemma
\ref{lm-normalspace}, a geodesic chart $\{ V_p, G_p, \pi_p \}$ can
be restricted to a chart $\{ C_p^+, G_p, \phi_p \}$ where $C_p^{+} =
\mathbf{B}_0^{n-1}(r/2) \times [0,\epsilon_p)$, $\phi_p$ is the
restriction of $\pi_p$ to $C_p^+$, and $\phi_p
\left(\mathbf{B}_0^{n-1} \times \{ 0 \} \right) =
\partial \phi_p(C_p^+)$.  We will refer to such a chart as a
\emph{boundary product chart} for $Q$.

It follows, in particular, that there is a neighborhood of $\partial
Q$ in $Q$ that is diffeomorphic to $\partial Q \times [0, \epsilon]$
for some $\epsilon > 0$ and that the metric respects the product
structure. This can be shown by forming a cover of $\partial Q$ of
sets uniformized by charts of the form $\{C^{+}, G_p, \psi_p\}$,
choosing a finite subcover, and setting $\epsilon = \mbox{min} \{
\epsilon_p/2 \}$.

\begin{lemma}
\label{lm-smoothifold}

The glued set $\hat{Q}$, i.e. the set of equivalence classes under
the identification made by Equation \ref{eq-glue}, may be made into
a smooth orbifold containing diffeomorphic copies of both $Q$ and
$Q^\prime$ such that $Q \cap Q^\prime = \partial Q = \partial
Q^\prime$.

\end{lemma}

{\it \noindent Proof:} For each point $p \in \partial Q$, form a
boundary product chart $\{ C_p^+, G_p, \phi_p \}$. Then glue each
chart of the boundary of $Q$ to its corresponding chart of
$Q^\prime$ in the following way. Let $\alpha: \mathbb{R}^n \to
\mathbb{R}^n$ be the reflection that sends $e_n \mapsto -e_n$ and
fixes all other coordinates.  A point $p$ in the boundary is
uniformized by two corresponding boundary product charts on either
side of $\partial Q$, $\{C_p^+, G_p, \phi_p \}$ and
$\{C_p^{+\prime}, G_p^\prime, \phi_p^\prime \}$. From these two
charts, a new chart $\{C_p, G_p, \psi_p \}$ for a neighborhood of
$p$ in $\hat{Q}$ is constructed where $C_p =
\mathbf{B}_{0}^{n-1}(r/2) \times (-\epsilon_p, \epsilon_p)$, and
\[
    \psi_p (x) = \left\{\begin{array}{ll}
    \phi_p (x),   & x_n \geq 0, \\
    \phi_p^\prime\circ\alpha (x),  & x_n < 0.
    \end{array}\right.
\]

These charts cover a neighborhood of $\partial Q = \partial
Q^\prime$ in $\hat{Q}$. By taking a geodesic chart at each point on
the interiors of $Q$ and $Q^\prime$ together with these new charts,
the entire set $\hat{Q}$ is covered.  Injections of charts at points
in the interior of $Q$ or $Q^\prime$ into charts of the form $\{
C_p^+, G_p, \phi_p \}$ induce injections into charts $\{ C_p, G_p,
\psi_p \}$.  Hence, $\hat{Q}$ is given the structure of a smooth
orbifold with the desired properties.

{\hfill $\square$}

Again, it follows that a neighborhood of $\partial Q \subset
\hat{Q}$ admits a tubular neighborhood diffeomorphic to $\partial Q
\times [-\epsilon, \epsilon]$ such that the metric respects this
product structure.


\section{The Morse Index of a Vector Field on an Orbifold}
\label{sec-morseindex}

The definitions of Morse Index and relation to the topological index
of a vector field extend readily to orbifolds, which we now
describe.

Let $Q$ be a compact orbifold with or without boundary, and let $X$
be a vector field on $Q$ that does not vanish on the boundary.  If
$X(p) = 0$ for $p \in Q$, then we say that $p$ is a
\emph{non-degenerate} zero of $X$ if there is a chart $\{ V, G, \pi
\}$ for a neighborhood $U_p$ of $p$ and an $x \in V$ with $\pi(x) =
p$ such that $\pi^\ast X$ has a non-degenerate zero at $x$; i.e.
$D(\pi^\ast X)_x$ has trivial kernel.  As in the manifold case,
non-degenerate zeros are isolated in charts and hence isolated on
$Q$. The Morse Index $\lambda (\pi^\ast X; x)$ of $\pi^\ast X$ at
$x$ is defined to be the number of negative eigenvalues of
$D(\pi^\ast X)_x$ (see \cite{Milnormorse}).  Since the Morse Index
is a diffeomorphism invariant, this index does not depend on the
choice of chart nor on the choice of $x$.  Since the
isomorphism-class of the isotropy group does not depend on the
choice of $x$, the expression $|G_p|$ is well-defined.  Hence, for
simplicity, we may restrict to charts of the form $\{ V_p, G_p,
\pi_p \}$ where $\pi_p(0) = p$ and $G_p$ acts linearly. We define
the \emph{orbifold Morse Index} of $X$ at $p$ to be
\[
    \lambda^{orb}(X; p)  =   \frac{1}{|G_p|} \lambda (\pi_p^\ast X; 0).
\]
Note that this index differs from that recently defined in
\cite{hepworth}; however, it is sufficient for our purposes. We have
\[
\begin{array}{rcl}
    \mathfrak{Ind}^{orb}(X;p)   &=&     \frac{1}{|G_p|} \mathfrak{Ind}(\pi^\ast X; 0)       \\\\
                                &=&     \frac{1}{|G_p|} (-1)^{\lambda(\pi^\ast X; 0)}.
\end{array}
\]

Suppose $X$ has only non-degenerate zeros on $Q$.  For each $\lambda
\in \{ 0, 1, \ldots , n\}$, we let $\{ p_i : i = 1, \ldots ,
k_\lambda \}$ denote the points in $Q$ at which the pullback of $X$
in a chart has Morse Index $\lambda$. Then we let
\[
    C_\lambda = \sum\limits_{i=1}^{k_\lambda} \frac{1}{|G_{p_i}|}
\]
count these points, where the orbifold-contribution of each zero
$p_i$ is $\frac{1}{|G_{p_i}|}$.  Note that as non-degenerate zeros
are isolated, there is a finite number on $Q$.

Then, as in the manifold case, if we define
\[
    \Sigma^{orb}(X; Q) =     \sum\limits_{\lambda = 0}^n (-1)^\lambda
     C_\lambda,
\]
we have
\[
\begin{array}{rcl}
    \Sigma^{orb}(X;Q)&=&     \sum\limits_{\lambda = 0}^n (-1)^\lambda \sum\limits_{i=1}^{k_\lambda} \frac{1}{|G_{p_i}|}
                                                                                     \\\\
                    &=&     \sum\limits_{p \in Q, X(p) = 0} \frac{1}{|G_p|} (-1)^{\lambda(\pi_p^\ast X; 0)}
                                                                                     \\\\
                    &=&     \sum\limits_{p \in Q, X(p) = 0}
                            \mathfrak{Ind}^{orb}(X;p)                                \\\\
                    &=&     \mathfrak{Ind}^{orb}(X;Q).
\end{array}
\]
In the case that $Q$ is closed, this quantity is equal to
$\chi_{orb}(Q)$ by Theorem \ref{thrm-satakeindex}.

We summarize these results as follows.

\begin{proposition}
\label{prop-morseindex}

Let $X$ be a smooth vector field on the compact orbifold $Q$ that
has non-degenerate zeros only, none of which occurring on $\partial
Q$. Then
\[
    \Sigma^{orb}(X; Q) = \mathfrak{Ind}^{orb}(X; Q).
\]
If $\partial Q = \emptyset$, then
\[
    \Sigma^{orb}(X; Q) = \chi_{orb}(Q).
\]

\end{proposition}

\begin{remark}
\label{rem-approximation}

If $Q$ is a compact orbifold (with or without boundary) and $X$ a
smooth vector field on $Q$ that is nonzero on some compact subset
$\Gamma$ of the interior of $Q$, then $X$ may be perturbed smoothly
outside of a neighborhood of $\Gamma$ so that it has only isolated,
nondegenerate zeros.  This is shown in \cite{wanerwu} for the case
of a smooth global quotient $M/G$ using local arguments, and so it
extends readily to the case of a general orbifold by working in
charts.

\end{remark}


\section{Proof of Theorem \ref{thrm-mainresult}}
\label{sec-mainresult}

Let $Y$ be a vector field in generic contact with $\partial Q$ that
has isolated zeroes on the interior of $Q$. Define $\hat{Y}$ on
$\hat{Q}$ by letting $\hat{Y}$ be $Y$ on each copy of $Q$.
Unfortunately, $\hat{Y}$ has conflicting definitions along the old
boundary $\partial Q$.  However, as in the manifold case treated in
\cite{Pugh}, the vector field may be perturbed near the boundary to
form a well-defined vector field using the product structure.  We
give an adaptation of Pugh's result to orbifolds.

\begin{proposition}
\label{prop-approxvfieldconstruction}

Given a smooth vector field $Y$ in generic contact with $\partial Q$
and with isolated zeros, none of which lie on $\partial Q$, there is
a smooth vector field $X$ on the double $\hat{Q}$ such that
\begin{itemize}
\item   Outside of a tubular neighborhood $P_\epsilon$ of $\partial Q$ containing none of the zeros of $Y$,
        $X$ coincides with $Y$ on $Q$ and $Q^\prime$;

\item   $X|_{\partial Q}$ is tangent to $\partial Q$,

\item   On $\Gamma^1$, $X$ coincides with $Y$ and in particular,
        defines the same $\Gamma^i$, $R_-^i$, and $R_+^i$ for $i > 1$; and

\item   The zeros of $X$ are those of $Y$ on the interior of $Q$ and
        $Q^\prime$ and a collection of isolated zeros on
        $\partial Q$ which are non-degenerate as zeros of $X|_{\partial Q}$.
\end{itemize}

\end{proposition}

{\it \noindent Proof:} As above, $\hat{Y}$ is defined everywhere on
$\hat{Q}$ except on the boundary.  Let $P_\epsilon$ be a normal
tubular $\epsilon$-neighborhood of $\partial Q$ in $\hat{Q}$ of the
form $\partial Q \times [-\epsilon, \epsilon]$ which we parameterize
as $\{ (x, v) : x \in \partial Q, v \in [-\epsilon, \epsilon] \}$.
We assume that $P_\epsilon$ is small enough so that it does not
contain any of the zeros of $\hat{Y}$. On $P_\epsilon$, decompose
$\hat{Y}$ respecting the product structure of $P_\epsilon$ into
\[
\hat{Y} = \hat{Y}_h + \hat{Y}_v
\]
These are the horizontal and vertical components of $\hat{Y}$,
respectively.  The horizontal component $\hat{Y}_h$ is well-defined,
continuous, and smooth when restricted to the boundary.  However,
$\hat{Y}_v$ has conflicting definitions on the boundary, although
they only differ by a sign. Note that the restriction of $\hat{Y}_h$
to $\partial Q$ may not have isolated zeros. However, as $Y$ does
not have zeros on $\partial Q$ and $\hat{Y}_h \equiv Y$ on
$\Gamma^1$, none of the zeros of $\hat{Y}_h|_{\partial Q}$ occur on
$\Gamma^1$.

Define $Z_h$ to be a smooth vector field on $\partial Q$ that
coincides with $\hat{Y}_h$ on an open subset of $\partial Q$
containing $\Gamma^1$ and has only non-degenerate zeros (see Remark
\ref{rem-approximation}).  Let $f(x, v)$ be the parallel transport
of $Z_h(x, 0)$ along the geodesic from $(x,0 )$ to $(x,v)$, and then
$Z_h$ is a horizontal vector field on $P_\epsilon$.  For $s \in (0,
\epsilon)$, let $\phi_s:\mathbb{R}\to [0, 1]$ be a smooth bump
function which is one on $[-s/2, s/2]$ and zero outside of $[s, s]$.

Define the vector field $X_s$ to be $\hat{Y}$ outside of
$P_\epsilon$ and
\[
    X_s(x, v)   =   \phi_s(v) \left(f(x,v) + |v|\hat{Y}_v(x,v)\right)
    + (1-\phi_s(v))\hat{Y}(x,v)
\]
on $P_\epsilon$.  Note that $X_s$ is smooth.  By picking $s$
sufficiently small, it may be ensured that the zeroes of $X$ are the
zeroes of $\hat{Y}$ and the zeroes of $Z_h|_{\partial{Q}}$ only. We
prove this as follows.

On points, $(x, v)$ where $x \in \Gamma^1$ and $|v| \leq s$, the
horizontal component of $X$ is $\phi_s(v) f(x, v) + (1 -
\phi_s(v))\hat{Y}_h(x, v)$.  Note that $f(x, 0) = \tilde{Y}_h(x, 0)$
for $x \in \Gamma^1$ and $f(x, 0) \neq 0$ on $\Gamma^1$.  Let $m
> 0$ be the minimum value of $\| f(x, 0) \|$ on the compact set
$\Gamma^1$, and then as $\Gamma^1 \times [-\epsilon, \epsilon]$ is
compact and $\tilde{Y}_h(x, v)$ continuous, there is an $s_0$ such
that
\[
        \| \hat{Y}_h(x, 0) - \hat{Y}_h(x, v) \|
        =   \| f(x, 0) - \hat{Y}_h(x, v) \| < m/2
\]
whenever $|v| < s_0$.  Hence, for such $v$ and for any $t \in [0,
1]$,
\[\
\begin{array}{rcl}
    \left\| t f(x, v) + (1 - t)\hat{Y}_h(x, v) \right\|
        &=&     \left\| \hat{Y}_h(x, v) + t [f(x, v) - \hat{Y}_h(x, v)] \right\|
                        \\\\
        &\geq&  \left \|\hat{Y}_h(x, v) \right\| - t \left\|f(x, v) - \hat{Y}_h(x, v) \right\|
                        \\\\
        &>&     m - \frac{tm}{2}
                        \\\\
        &\geq&  \frac{m}{2} \; >   \;  0.
\end{array}
\]
Therefore, the horizontal component is nonvanishing, implying that
$X_s(v, h)$ does not vanish here.

Now let $\{ x_i : i = 1, \ldots, k \}$ be the zeros of $Z_h$ on
$\partial Q$.  Each $x_i$ is contained in a ball $B_{\epsilon_i}
\subset \partial Q$ whose closure does not intersect $\Gamma^1$.
Hence, $\hat{Y}_v(x, 0) \neq 0$ on each $B_{\epsilon_i}$. Therefore,
for each $i$, there is an $s_i$ such that $\hat{Y}_v(x, v) \neq 0$
on $B_{\epsilon_i} \times (-s_i, s_i)$.  This implies that the
vertical component of $X_s(x, v)$, and hence $X_s(x,v)$ itself, does
not vanish on $B_{\epsilon_i} \times (-s_i, s_i)$ except where $v =
0$; i.e. on $\partial Q$.

Letting $s$ be less than the minimum of $\{ s_0, s_1, \ldots,
s_k\}$, we see that $X_s$ does not vanish on $P_\epsilon$ except on
$\partial Q$, where it coincides with $Z_h$. Therefore, $X = X_s$ is
the vector field which was to be constructed.

{\hfill $\square$}

It follows that the index of the vector field $X$ constructed in the
proof of Proposition \ref{prop-approxvfieldconstruction} is
\begin{equation}
\label{eq-indexstep1}
    \mathfrak{Ind}^{orb}(X;Q)   =   2\mathfrak{Ind}^{orb}(Y;Q)
            + \sum_{p \in \partial Q}\mathfrak{Ind}^{orb}(X;p)
\end{equation}
Let $p$ be a zero of $X$ on $\partial{Q}$, i.e. it is a zero of
$Z_h$.  We will write the index of $X$ at $p$ in terms of the index
of $Z_h$.

Because of Lemma \ref{lm-normalspace}, the isotropy group of $p$ as
an element of $Q$ is the same as the isotropy group of $p$ as an
element of $\partial{Q}$, and so we may refer to $G_p$ without
ambiguity.  About a neighborhood of $p$ in $Q$ small enough to
contain no other zeros of $X$, choose a boundary product chart
$\{C_p^+, G_p, \phi_p \}$. Then, as in Lemma \ref{lm-smoothifold},
$\{ C, G_p, \psi_p \}$ forms a chart about $p$ in $\hat{Q}$.  The
product structure $(y, w)$ of $C_p = \mathbf{B}_0^{n - 1}(r/2)
\times (-\epsilon_p, \epsilon_p)$ coincides with that of
$P_\epsilon$ near the boundary, so within the preimage of $\partial
Q \times [-s/2, s/2]$ by $\psi_p$, we have that
\[
    \psi_p^\ast X = \psi_p^\ast f + |w| \psi_p^\ast \hat{Y}_v.
\]
Note that $\psi_p(0,0) = (p)$, and then
\[
\begin{array}{rcl}
    D(\psi_p^\ast X)_{(0, 0)}  &=& \left(  \begin{array}{cc}
                    D(\psi_p^\ast Z_h)_0
                &   \left(\frac{\partial \psi_p^\ast f}{\partial w} \right)_0
                                \\\\
                    D\left((|w| \psi_p^\ast \hat{Y}_v)|_{\partial C_p }\right)_0
                &   \left(\frac{\partial}{\partial w} |w| \psi_p^\ast \hat{Y}_v \right)_0
                \end{array} \right)
        \\\\
            &=& \left(  \begin{array}{cc}
                    D(\psi_p^\ast Z_h)_0
                &   0           \\\\
                    0   &   \psi_p^\ast \hat{Y}_v(0, 0)
                \end{array} \right) .
\end{array}
\]
As $\psi_p^\ast \hat{Y}_v(0, 0)$ is positive if $p \in R_+^1$ and
negative if $p \in R_-^1$, it is seen that
\[
    \lambda\left(\psi_p^\ast X ; (0,0) \right) =
    \left\{\begin{array}{ll}
        \lambda\left(\psi_p^\ast X|_{\partial C_p }; 0 \right),
                &       p \in R_{+},            \\
        \lambda\left(\psi_p^\ast X|_{\partial C_p }; 0 \right) + 1,
                &       p \in R_{-} .
\end{array}\right.
\]
Hence
\[
    \mathfrak{Ind}\left(\psi_p^\ast X ;(0,0) \right) =
    \left\{\begin{array}{ll}
        \mathfrak{Ind}\left(\psi_p^\ast Z_h|_{\partial C_p }; 0\right),
                &       p \in R_{+},            \\
        -\mathfrak{Ind}\left(\psi_p^\ast Z_h|_{\partial C_p }; 0 \right),
                &       p \in R_{-} .
\end{array}\right.
\]
Therefore, for $p \in R_+$,
\[
\begin{array}{rcl}
    \mathfrak{Ind}^{orb}(X, p)
        &=&     \frac{1}{|G_p|} \mathfrak{Ind}\left(\psi_p^\ast X ;0 \right)
                                                                                \\\\
        &=&     \frac{1}{|G_p|} \mathfrak{Ind}\left(\psi_p^\ast Z_h|_{\partial C^+}; 0 \right)
                                                                                \\\\
        &=&     \mathfrak{Ind}^{orb}\left(Z_h ;p \right)
\end{array}
\]
and similarly
\[
    \mathfrak{Ind}^{orb}\left(X ;p \right)
    = - \mathfrak{Ind}^{orb}\left(Z_h ;p \right)
\]
for $p \in R_-$.

With this, Equation \ref{eq-indexstep1} becomes
\[
    \mathfrak{Ind}^{orb}(X; \hat{Q})
    =   2\mathfrak{Ind}^{orb}(Y;Q)
    +   \mathfrak{Ind}^{orb}(Z_h; R_+)
    -   \mathfrak{Ind}^{orb}(Z_h; R_-).
\]

By Theorem \ref{thrm-satakeindex} and Equation
\ref{eq-additiveeulerchar}, $\mathfrak{Ind}^{orb}(X; \hat{Q}) =
2\chi_{orb}(Q) - \chi_{orb}(\partial Q)$, so that
\[
    2\chi_{orb}(Q) - \chi_{orb}(\partial Q)
    =
    2\mathfrak{Ind}^{orb}(Y;Q)
    +   \mathfrak{Ind}^{orb}(Z_h; R_+)
    -   \mathfrak{Ind}^{orb}(Z_h; R_-).
\]
Note that $\partial Q$ is also a closed orbifold, so
\[
\begin{array}{rcl}
    \chi_{orb}(\partial Q)
        &=&         \mathfrak{Ind}^{orb}(X; \partial Q)     \\\\
        &=&         \mathfrak{Ind}^{orb}(X; R_{+}) - \mathfrak{Ind}^{orb}(X;
                    R_{-}).
\end{array}
\]
Hence,
\begin{equation}
\label{eq-indexstep2}
\begin{array}{rcl}
\mathfrak{Ind}^{orb}(Y;Q)
            &=& \chi_{orb}(Q) + \frac{1}{2}\left(-\chi_{orb}(\partial Q)  + \mathfrak{Ind}^{orb}(X; R_{-})
                - \mathfrak{Ind}^{orb}(X; R_{+}) \right)
                \\\\
            &=& \chi_{orb}(Q) + \frac{1}{2}\left[-\chi_{orb}(\partial Q)  + 2\mathfrak{Ind}^{orb}(X; R_{-}) \right.
                \\\\
            &&  \left.
                - \left(\mathfrak{Ind}^{orb}(X; R_{+}) + \mathfrak{Ind}^{orb}(X; R_{-}) \right)\right]
                \\\\
            &=& \chi_{orb}(Q) + \frac{1}{2}\left(-2\chi_{orb}(\partial Q)  + 2\mathfrak{Ind}^{orb}(X;R_{-})\right)
                \\\\
            &=& \chi_{orb}(Q) - \chi_{orb}(\partial Q)  + \mathfrak{Ind}^{orb}(X; R_{-})
                \\\\
            &=& \chi_{orb}(Q, \partial Q) + \mathfrak{Ind}^{orb}(X; R_{-})
\end{array}
\end{equation}

Because $X$ coincides with $Y$ on $\Gamma^1$, it defines the same
$\Gamma^{i}$ that $Y$ does.  Since $X$ is a smooth vector field
defined on $R_{-}^1$ that does not vanish on $\partial R_-^1 =
\Gamma^1$, we may recursively apply this formula to higher and
higher orders of $R^{i}_{-}$ until $R^{i}_{-}$ is empty, and there
will no longer be an index sum term. Hence,
\[
    \mathfrak{Ind}^{orb}(X; R_{-}) = \sum_{i = 1}^n \chi_{orb}(R_{-}^i,
    \Gamma^i).
\]
Along with Equation \ref{eq-indexstep2}, this completes the proof of
Theorem \ref{thrm-mainresult}.

{\hfill $\square$}

Let $\tilde{Q}$ denote the inertia orbifold of $Q$ and $\pi :
\tilde{Q} \rightarrow Q$ the projection (see
\cite{chenruanorbcohom}). It is shown in \cite{seaton1} that a
vector field $Y$ on $Q$ induces a vector field $\tilde{Y}$ on
$\tilde{Q}$, and that $\tilde{Y}(p, (g)) = 0$ if and only if $Y(p) =
0$.

For each point $p \in Q$ and $g \in G_p$, a chart $\{ V_p, G_p,
\pi_p\}$ induces a chart $\{ V_p^g, C(g), \pi_{p, g} \}$ at $(p,
(g)) \in \tilde{Q}$ where $V_p^g$ denotes the points in $V_p$ fixed
by $g$ and $C(g)$ denotes the centralizer of $g$ in $G_p$. Clearly,
$\partial V_p^g = (\partial V_p) \cap V_p^g$.  An atlas for
$\tilde{Q}$ can be taken consisting of charts of this form, so it is
clear that $\partial \tilde{Q} = \widetilde{\partial Q}$.

Let $p \in \partial Q$ and pick a boundary product chart $\{ C_p^+,
G_p, \phi_p \}$.  Then for $g \in G_p$, there is a chart $\{
(C_p^+)^g, C(g), \phi_{p, g} \}$ for $(p, (g)) \in \tilde{Q}$.  As
the normal component to the boundary of $C_p^+$ is $G_p$-invariant,
\[
\begin{array}{rcl}
    (C_p^+)^g       &=&     \left(\mathbf{B}_0^{n-1}(r/2) \times [0, \epsilon_p)
                            \right)^g                       \\\\
                    &=&     \left(\mathbf{B}_0^{n-1}(r/2) \right)^g \times [0,\epsilon_p),
\end{array}
\]
and so
\[
    T_0 (C_p^+)^g       =       T_0 \left(\mathbf{B}_0^{n-1}(r/2)
                                \right)^g \times \R.
\]
It follows that $\tilde{Y}$ points out of $\partial \tilde{Q}$ at
$(p, (g))$ if and only if $Y$ points out of $\partial Q$ at $p$.
With this, applying Theorem \ref{thrm-mainresult} to $\tilde{Y}$
yields
\begin{equation}
\label{eq-inert}
\begin{array}{rcl}
    \mathfrak{Ind}^{orb}(\tilde{Y}; \tilde{Q})
        &=& \chi_{orb}(\tilde{Q}, \partial \tilde{Q}) +
            \sum\limits_{i=1}^n \chi_{orb}\left(\widetilde{R_-^i},
            \widetilde{\Gamma^i}\right)         \\\\
        &=& \chi_{orb}(\tilde{Q}) - \chi_{orb}(\partial \tilde{Q}) +
            \sum\limits_{i=1}^n
        \chi_{orb}\left(\widetilde{R_-^i}\right) - \chi_{orb}\left(
            \widetilde{\Gamma^i}\right).
\end{array}
\end{equation}
Each of the $\Gamma^i$ and $\partial Q$ are closed orbifolds, so it
follows from the proof of Theorem 3.2 in \cite{seaton1} (note that
the assumption of orientability is not used to establish this
result) that
\[
    \chi_{orb}(\widetilde{\Gamma^i}) = \chi(\mathbb{X}_{\Gamma^i})
\]
and
\begin{equation}
\label{eq-orbtopboundary}
    \chi_{orb}(\partial \tilde{Q}) = \chi_{orb}(\widetilde{\partial Q}) = \chi(\mathbb{X}_{\partial Q})
\end{equation}
where $\mathbb{X}_{\Gamma^i}$ and $\mathbb{X}_{\partial Q}$ denote
the underlying topological spaces of $\Gamma^i$ and $\partial Q$,
respectively, and $\chi$ the usual Euler characteristic.

Letting $\hat{Q}$ denote, as above, the double of $Q$, it is easy to
see that $\hat{\tilde{Q}} = \tilde{\hat{Q}}$.  Hence, applying the
same result to $\hat{\tilde{Q}}$ yields
\begin{equation}
\label{eq-doubleinert}
\begin{array}{rcl}
    \chi(\mathbb{X}_{\hat{Q}})
        &=&     \chi_{orb}\left( \tilde{\hat{Q}} \right)            \\\\
        &=&     \chi_{orb} \left( \hat{\tilde{Q}} \right)           \\\\
        &=&     2\chi_{orb}(\tilde{Q}) - \chi_{orb}(\partial \tilde{Q})
\end{array}
\end{equation}
However, as
\[
\begin{array}{rcl}
    \chi(\mathbb{X}_{\hat{Q}})      &=&     2\chi(\mathbb{X}_{Q}) -
                                            \chi(\mathbb{X}_{\partial Q})
            \\\\
      &=&     2\chi(\mathbb{X}_{Q}) - \chi_{orb}(\partial \tilde{Q})
\end{array}
\]
it follows from Equation \ref{eq-doubleinert} that
$\chi_{orb}(\tilde{Q}) = \chi(\mathbb{X}_{Q})$. The same holds for
each $R_-^i$ so that Equation \ref{eq-inert} becomes the following.

\begin{corollary}
\label{cor-inertiaversion}

Let $Q$ be an $n$-dimensional smooth, compact orbifold with
boundary, and let $Y$ be a smooth vector field on $Q$.  If
$\tilde{Y}$ denotes the induced vector field on $\tilde{Q}$, then
\[
    \mathfrak{Ind}^{orb}(\tilde{Y}; \tilde{Q})
        = \chi(\mathbb{X}_Q, \mathbb{X}_{\partial Q}) +
        \sum\limits_{i=1}^n \chi(\mathbb{X}_{R_-^i},
        \mathbb{X}_{\Gamma^i}).
\]

\end{corollary}


\bibliographystyle{amsplain}

\end{document}